\documentclass{IEEEtran}
\usepackage{cite}
\usepackage{amsmath,amssymb,amsfonts,dsfont}
\usepackage{algorithmic}
\usepackage{graphicx}
\usepackage{textcomp}
\usepackage{algorithm}
\usepackage{subfigure}
\newcommand{\tabincell}[2]{\begin{tabular}{@{}#1@{}}#2\end{tabular}}
\usepackage{soul}
\usepackage{xcolor}
\usepackage{epstopdf}
\newtheorem{thm}{Theorem}
\newtheorem{remark}{Remark}
\newtheorem{prop}{Proposition}
\newtheorem{coro}{Corollary}
\newtheorem{assumption}{Assumption}
\makeatletter
\def\BibTeX{{\rm B\kern-.05em{\sc i\kern-.025em b}\kern-.08em
    T\kern-.1667em\lower.7ex\hbox{E}\kern-.125emX}}
\begin{document}

\title{Efficient Sampling Policy for Selecting a Good Enough Subset}

\author{Gongbo Zhang, Bin Chen, Qing-shan Jia and Yijie Peng
\thanks{This work was supported in part by the National Natural Science Foundation of China (NSFC) under Grants 62073182, 71901003 and 72022001.}
\thanks{Gongbo Zhang and Yijie Peng are with Department of Management Science and Information Systems, Guanghua School of Management, Peking University, Beijing 100084, China (e-mail: gongbozhang@pku.edu.cn; pengyijie@pku.edu.cn). }
\thanks{Bin Chen is with College of Systems Engineering, National University of Defense Technology, Changsha 410073, China (e-mail: nudtcb9372@gmail.com).}
\thanks{Qing-Shan Jia is with the Center for Intelligent and Networked Systems (CFINS), Department of Automation, Beijing National Research Center for Information Science and Technology (BNRist), Tsinghua University, Beijing 100084, China (e-mail: jiaqs@tsinghua.edu.cn).}
\thanks{Corresponding authors: Bin Chen, Yijie Peng}
}

\maketitle

\begin{abstract}
The note studies the problem of selecting a good enough subset  out of a finite number of alternatives under a fixed simulation budget. Our work aims to maximize the posterior probability of correctly selecting a good subset. We formulate the dynamic sampling decision as a stochastic control problem in a Bayesian setting. In an approximate dynamic programming paradigm, we propose a sequential sampling policy based on value function approximation. We analyze the asymptotic property of the proposed sampling policy. Numerical experiments demonstrate the efficiency of the proposed procedure.
\end{abstract}


\section{Introduction}

Simulation optimization (SO) has been widely used in analyzing modern complex system, e.g., healthcare, transportation, manufacturing and supply chain systems \cite{fu2015handbook}. Simulation experiments are often expensive and time-consuming since a large number of simulation replications are required in order to achieve an accurate estimate of the system performance \cite{chen2011stochastic}. In this note, we consider the problem of selecting a subset of size $m$ containing the best of $k$ alternatives, $1 \le m < k$, $m,k \in \mathbb{Z}^+$. We call such subset a good enough subset in this work. The true performance of each alternative is unknown in practice, which may be estimated by Monte Carlo simulation under a finite sampling budget $T$.

{When $m = 1$, the problem reduces to the ranking and selection (R\&S) for finding the best.} R\&S has been widely studied in simulation, which aims to allocate simulation replications to alternatives for efficiently selecting the best alternative. {There are fixed-precision and fixed-budget procedures for R\&S~\cite{hong2020review}, where the former allocates samples to guarantee a probability of correct selection (PCS) up to a certain level, whereas the latter optimizes performance metrics under a fixed simulation budget constraint. A well-researched paradigm is an indifference zone (IZ) framework \cite{rinott1978two,kim2001fully}, which aims to guarantee a PCS level in the least favorable configuration. The sampling procedures in the IZ framework tend to allocate more replications than necessary in practice. Well-known methods to enhance the efficiency for finding the best alternative include} optimal computing budget allocation (OCBA) \cite{chen2000simulation}, expected value of information (EVI) \cite{chick2001new,chick2010sequential}, knowledge gradient (KG) \cite{gupta1996bayesian,frazier2009knowledge}, expected improvement (EI) \cite{jones1998efficient,ryzhov2016convergence} and asymptotically optimal allocation procedure (AOAP) \cite{peng2018ranking}. OCBA was initially a two-stage procedure, and becomes fully sequential by combining with a certain sequential rule, e,g,. ``most starving rule"~\cite{chen2011stochastic}. The fully sequential procedures usually lead to higher PCS under the specified simulation budget, or achieve the same level of PCS using fewer simulation replications than the two-stages procedures~\cite{chen2000simulation}. AOAP is derived under a stochastic control framework, and it is proved to achieve the asymptotically optimal sampling ratios for selecting the best with a normal sampling distribution \cite{glynn2004}. {Our problem softens the optimization problem by allowing a subset containing the best to be acceptable}, so that decision can be made in a more flexible way. Selecting a good enough subset instead of a single alternative as the best can dramatically reduce the computational cost and this problem formulation is particularly useful for many engineering applications when finding the best is not necessary and simulation is expensive. {For example, in an emergency department healthcare staffing problem, the number of staffing designs is large and circumstances of patients are complex, therefore, finding a good enough subset of designs can reduce the time for decision-making. In addition, good enough subset selection can efficiently screen out some alternatives at the first stage before seeking an optimal decision accurately at the second stage.}

{Compared with the R\&S for selecting the best alternative, the literature of good enough subset selection problem is still sparse. \cite{koenig1985procedure} develops a two-stage indifference zone procedure to select a subset of size $m$ containing the $\ell$ best of $k$ alternatives, $1 \le \ell \le m < k$. Following the pioneer work \cite{koenig1985procedure} using a least favourable configuration, \cite{chen2009subset} proposes a two-stage procedure to select a subset of size $m$ containing at least $c$ of the $\ell$ best of $k$ alternatives with unknown means and unknown variances. The most relevant work to us is \cite{gao2017new}, which studies the problem of selecting $\ell$ good enough alternatives out of $m$ acceptable alternatives from $k$ alternatives under a fixed simulation budget constraint. Yet strictly speaking, the objective in \cite{gao2017new} is to correctly select a subset of size $m$ containing the $\ell$ best of $k$ alternatives. By relaxing the objective function using a Bonferroni inequality and asymptotic analysis, \cite{gao2017new} proposes a sampling procedure that screens out  some alternatives at the first stage. %and neither finite nor asymptotic performances of the proposed method are analyzed.
In contrast to \cite{gao2017new}, we formulate the sampling decision as a stochastic dynamic programming problem. Although we do not screen out any alternative, it can be proved that following our procedure, sampling ratios of $\left(m-1\right)$ alternatives converge to 0 almost surely as the simulation budget goes to infinity, which explains why selecting a good enough subset could reduce the computational cost.}

{Our work is related to the literature on subset selection, which selects a subset of random size that contains the best alternative with a high probability. Different than our work, subset selection procedures guarantee the PCS level and lead to a small expected subset size \cite{gupta1965some,boesel2003using,eckman2020revisiting}, or they minimize various loss functions under a Bayesian framework \cite{studden1967selecting,deely1968properties,mieseke1979bayesian}. Another line of literature is selection of good alternative, which finds an alternative whose performance is within a given tolerance range from that of the best alternative \cite{fan2016indifference,eckman2018guarantees}, or an alternative whose performance is within the top-$m\%$ or top-$m$ alternatives\cite{ho2000ordinal}.% Apropos of that line of literature, our problem enables the best alternative in the good enough subset. We draw connections and provide insights into application of our framework to the good alternative selection problem.
}

We formulate the  dynamic sampling decision as a stochastic  control problem in a Bayesian setting. With an approximation of the PCS, an efficient sequential sampling procedure is derived in an approximate dynamic programming (ADP) paradigm by maximizing a value function approximation (VFA) one-step look ahead. {The proposed asymptotically optimal allocation procedure for selecting a good enough subset (AOA-gs) is proved to be consistent, that is, the best alternative will eventually belong to the good enough subset almost surely as the number of simulation budget goes to infinity, and the asymptotic sampling ratios of AOA-gs are derived.} We test the performance of the proposed sequential AOA-gs procedure via a series of numerical experiments including an $\left(s,S\right)$-type inventory problem and an emergency department healthcare staffing problem. Numerical results show that AOA-gs outperforms the existing methods in all tested experiments. The advantage of AOA-gs appears to be more pronounced when the number of competing alternatives is large.

The rest of the note is organized as follows. Section \ref{sec2} formulates the problem for selecting a good enough alternative based on stochastic control framework. The dynamic sampling allocation scheme is proposed in \ref{sec3}, and Section \ref{sec4} provides numerical experiments. The last section concludes the note.

\section{Problem Formulation}\label{sec2}

Suppose there are $k$ alternatives and the performance of each alternative is measured by unknown mean ${\mu _i} \in \mathbb{R}$, $i=1,\cdots,k$, where $\mu_i$ is estimated by Monte Carlo simulation. We assume that $\mu_i \ne \mu_j$ for $i \ne j$, $i,j=1,\cdots,k$, which ensures that the performance of each alternative is distinguishable. The objective of the problem is to select a good enough subset of size $m$ containing the best alternative $\left\langle 1 \right\rangle$, where $m$ is specified and ${\left\langle i \right\rangle }$, $i=1,\cdots,k$ are indices ranked by system performances such that ${\mu _{\left\langle 1 \right\rangle }}  >  \cdots  > {\mu _{\left\langle k \right\rangle }}$. Let $T < \infty$ be the number of total simulation replications, and $X_{i,t}$, $0 < t \le T$ be the $t$-th independent and identically distributed (i.i.d.) replication for alternative $i$. Suppose that $X_t \mathop  = \limits^\Delta \left(X_{1,t},\cdots,X_{k,t}\right)$ follows a joint sampling distribution $Q\left( { \cdot ;\theta } \right)$, where $\theta  \in \Theta$ is a vector comprising of all unknown parameters in the parametric family, and in particular, {$\left\{\mu_1,\cdots,\mu_k\right\} \subset \theta$}. The replications across different alternatives are assumed to be independent. The unknown parameter $\theta$ follows a prior distribution $F\left( { \cdot ;{\zeta _0}} \right)$, where $\zeta_0$ contains all hyper-parameters for the parametric family of the prior distribution. In practice, some information gathered from previous studies, export opinions and individual studies can be incorporated into the prior.

The dynamic allocation policy is a sequence of mappings ${\mathcal{A}_T} \mathop  = \limits^\Delta \left( {{A_1}\left(  \cdot  \right), \cdots ,{A_T}\left(  \cdot  \right)} \right)$, {${A_t}:{\mathcal{E}_{t - 1}} \to \left\{ {1, \cdots ,k} \right\}$}, where $A_t\left(\cdot\right) \in \left\{1,\cdots,k\right\}$ allocates the $t$-th replication to an alternative based on the available information set $\mathcal{E}_{t-1}$ collected throughout the first $\left(t-1\right)$ allocated replications. Denote the sample observations of each alternative throughout $t$ steps as $X_i^{\left( t \right)} \mathop  = \limits^\Delta \left( {{X_{i,1}}, \cdots ,{X_{i,{t_i}}}} \right)$, $i=1,\cdots,k$, where ${t_i} = \sum\nolimits_{\ell = 1}^t {{A_{i,\ell}}} \left( {{\mathcal{E} _{\ell - 1}}} \right)$, ${A_{i,\ell}}\left( {{\mathcal{E} _{\ell - 1}}} \right) = \mathds{1}\left( {{A_\ell}\left( {{\mathcal{E} _{\ell - 1}}} \right) = i} \right)$ and $\mathds{1} \left(\cdot\right)$ is an indicator function that equals 1 if the event in the bracket is true. ${\mathcal{E} _t} = \{ {{\zeta _0},X_1^{\left( t \right)}, \cdots ,X_k^{\left( t \right)}} \}$ comprises all sample observations and prior information $\zeta_0$. {We determine the good enough subset after allocating $T$ replications by the estimated set ${\widehat{\mathcal{F}}_T^m}\mathop  = \limits^\Delta  \left\{ {{\left\langle 1 \right\rangle}_T,\cdots ,{\left\langle m \right\rangle}_T } \right\}$, where ${\left\langle i \right\rangle}_T$, $i=1,\cdots,k$ are indices ranked by posterior means. A correct selection occurs when ${\left\langle 1 \right\rangle} \in {\widehat{\mathcal{F}}_T^m}$.
\begin{prop} The posterior PCS for a good enough subset can be expressed as
\begin{equation}
\begin{aligned}\label{PCS}
{\rm {PCS}}_T & = \Pr \left\{\left.{\left\langle 1 \right\rangle} \in {\widehat{\mathcal{F}}_T^m} \right|{\mathcal{E}_T}\right\} \\
& = \Pr \bigg\{{\bigcup\nolimits_{i = 1}^m{\bigcap\nolimits_{j = m+1}^k {{\mu _{{{\left\langle {{i}} \right\rangle }_T}}} > {\mu _{{\left\langle j \right\rangle }_T}}} }} \bigg|{\mathcal{E}_T} \bigg\}.
\end{aligned}
\end{equation}
\end{prop}
The proof of the proposition can be found in the online appendix \cite{zhang2022online}.} The closed-form expression of (\ref{PCS}) is unknown, and we aim to find a dynamic sampling decision such that (\ref{PCS}) is maximized. Note that a lower bound of (\ref{PCS}) is
\begin{align}\label{lb}
\mathop {\max }\nolimits_{i = 1, \cdots ,m} {\Pr \left\{ {\left. \bigcap\nolimits_{j = m + 1}^k {{\mu _{{{\left\langle {i} \right\rangle }_T}}} > {\mu _{{{\left\langle j \right\rangle }_T}}}} \right|{\mathcal{E} _{T}} } \right\}}~,
\end{align}
which is an estimate of the posterior PCS (\ref{PCS}), and (\ref{PCS}) goes to one when (\ref{lb}) goes to one. In this work, we propose a dynamic allocation policy $\mathcal{A}_{T}$ to maximize (\ref{lb}). The dynamic sampling decision can be captured by a stochastic control problem. Under the Bayesian setting, we recursively define the expected payoff for $\mathcal{A}_T$ by ${V_T}\left( {{\mathcal{E} _T};{\mathcal{A}_T}} \right) \mathop  = \limits^\Delta \mathop {\max }\nolimits_{i = 1, \cdots ,m} \Pr \{ { \bigcap\nolimits_{j = m + 1}^k {{\mu _{{{\left\langle {i} \right\rangle }_T}}} > {\mu _{{{\left\langle j \right\rangle }_T}}}} |{\mathcal{E} _{T}} } \}$, and for $0 \le t < T$,
$${V_t}\left( {{\mathcal{E} _t};{\mathcal{A}_T}} \right) \mathop = \limits^\Delta {\left. {\mathds{E}\left[ {\left. {{V_{t + 1}}\left( {{\mathcal{E} _t} \cup \left\{ {{X_{i,t_i + 1}}} \right\};{\mathcal{A}_T}} \right)} \right|{\mathcal{E} _t}}\right]} \right|_{i = {A_{t + 1}}\left( {{\mathcal{E} _t}} \right)}}~.$$

Then the optimal allocation policy $\mathcal{A}_T^*$ can be well defined by ${\mathcal{A}_T^*}\mathop = \limits^\Delta  \arg \mathop {\max }\nolimits_{{\mathcal{A}_T}} {V_0}\left( {{\zeta _0};{\mathcal{A}_T}} \right)$. The stochastic control problem can be viewed as a Markov decision process (MDP) with $\left(T+1\right)$ stages, where $\zeta_0$ is the state at stage 0, and $\mathcal{E}_t$ is the state at stage $t$, $0 < t \le T$. The actions correspond to $A_{t+1}$ for $0 \le t < T$ and the terminal selection action for $t=T$. The transition is ${\mathcal{E}_t} \to {\mathcal{E}_{t+1}}$, where ${\mathcal{E}_{t+1}} \mathop  = \limits^\Delta  \{{\mathcal{E} _t} \cup \left\{{X_{i,{t_i + 1}}}\right\}\}$. The only nonzero reward is the terminal reward $V_T\left(\mathcal{E}_T\right)$. By induction, the solution to the stochastic control problem and the optimal policy for the Bellman equation of MDP are equivalent~\cite{bertsekas1995dynamic}. We can recursively solve the $\mathcal{A}_T^{*}$ by the Bellman equation: ${V_t}\left( {{\mathcal{E} _t}} \right) \mathop  = \limits^\Delta {\left. {\mathbb{E}\left[ {\left. {{V_{t + 1}}\left( {\mathcal{E} _{t+1}} \right)} \right|{\mathcal{E} _t}} \right]} \right|_{i = A_{t + 1}^*\left( {{\mathcal{E} _t}} \right)}}$, $0 \le t < T$, where $A_{t + 1}^*\left( {{\mathcal{E} _t}} \right) = \arg \mathop {\max }\nolimits_{i = 1, \cdots ,k} {\mathbb{E}\left[ {\left. {{V_{t + 1}}\left( {\mathcal{E} _{t+1}} \right)} \right|{\mathcal{E} _t}} \right]}$. Finding $\mathcal{A}_T^{*}$ through backward induction suffers from the curse-of-dimensionality~\cite{peng2018ranking}. To address the computational difficulty, we adopt an approximate dynamic programming (ADP) paradigm~\cite{powell2007approximate}, which makes dynamic decision based on VFA and keeps learning the value function with decisions moving forward.

Suppose ${X_{i,t}}\mathop  \sim \limits^{i.i.d.} N\left( {{\mu _i},\sigma _i^2} \right)$, $i=1,\cdots,k$, with unknown means and known variances. To make the proposed allocation policy amenable to practical implementation, we focus on the known variance
case and use the sample estimate as a plug-in for the true value. If the normal assumption is not satisfied, a macro replication obtained from a batched mean follows an approximately normal distribution by the central limit theorem. From~\cite{chen1996}, we have the following assumption:
\begin{assumption}
$\mu_i$, $i \in \left\{1,\cdots,k\right\}$, follows a conjugate normal prior distribution.
\end{assumption}

By \cite{degroot2005optimal}, the conjugate prior for a normal sampling distribution $N\left( {{\mu _i},\sigma _i^2} \right)$ with unknown mean and known variance is also a normal distribution $N ( {\mu _i^{\left( 0 \right)},{( {\sigma _i^{( 0 )}} )^2}} )$. The posterior distribution of $\mu_i$ is $N( {\mu _i^{\left( t \right)},{( {\sigma _i^{( t )}} )^2}} )$ , where
$$( \sigma _i^{\left( t \right)} )^2 = {\left( {\frac{1}{{{( {\sigma _i^{( 0 )}} )^2}}} + \frac{{{t_i}}}{\sigma_i^2}} \right)^{ - 1}}~,$$
$$\mu_i^{\left( t \right)} = {( \sigma _i^{\left( t \right)} )^2}\left( {\frac{{\mu _i^{\left( 0 \right)}}}{{{( {\sigma _i^{( 0 )}} )^2}}} + \frac{{{t_i}m_i^{\left( t \right)}}}{\sigma_i^2}} \right),\;m_i^{\left( t \right)} = \frac{\sum\nolimits_{\ell = 1}^{{t_i}} {{X_{i,\ell}}} } {{t_i}}~.$$

If $\sigma _i^{\left( 0 \right)} \to \infty $, $\mu_i^{\left(t\right)} = m_i^{\left(t\right)}$, and such a prior is uninformative. For a normal distribution with unknown variance, there is a normal-gamma conjugate prior~\cite{degroot2005optimal}. By conjugacy, $\mathcal{E}_t$ is completely determined by the posterior hyper-parameters, and the dimension of $\mathcal{E}_t$ is fixed at any step.

\section{Dynamic allocation policy}\label{sec3}

In this section, in order to derive a dynamic allocation procedure with an analytical formt, we adopt the VFA technique in \cite{peng2018ranking}, which uses a single feature of
the value function one-step ahead. Specifically, suppose any step $t$ could be the last step. The joint distribution of vector
$$\left( {{\mu _{{{\left\langle {i} \right\rangle }_t}}} - {\mu _{{{\left\langle {m + 1} \right\rangle }_t}}}}, \cdots ,{{\mu _{{{\left\langle {i} \right\rangle }_t}}} - {\mu _{{{\left\langle k \right\rangle }_t}}}}\right),\quad i\in\left\{1,\cdots,m\right\}~,$$
follows a joint normal distribution with mean vector
$$\left( {\mu _{_{{{\left\langle {i} \right\rangle }_t}}}^{\left( t \right)} - \mu _{{{\left\langle {m + 1} \right\rangle }_t}}^{\left( t \right)},  \cdots ,\mu _{_{{{\left\langle {i} \right\rangle }_t}}}^{\left( t \right)} - \mu _{{{\left\langle {k} \right\rangle }_t}}^{\left( t \right)}}\right),\quad i\in\left\{1,\cdots,m\right\}~,$$
and covariance matrix $\Sigma_i ={\widetilde{\Gamma}} '{\Lambda_i} {\widetilde{\Gamma}}$, where $^\prime$ denotes the transpose operation of the matrix, ${\widetilde{\Gamma}} \in \mathbb{R}^{\left( {k - m + 1} \right) \times \left( {k - m} \right)}$,
$${\widetilde{\Gamma}} \mathop  = \limits^\Delta {\left( {\begin{array}{*{20}{c}}
1&1&1& \cdots &1\\
{ - 1}&0&0& \cdots &0\\
0&{ - 1}&0& \cdots &0\\
 \vdots & \vdots & \vdots & \cdots & \vdots \\
0&0&0& \cdots &{ - 1}
\end{array}} \right)}~,$$
and $\Lambda_i \in \mathbb{R}^{\left(k-m+1\right) \times \left(k-m+1\right)}$ is a diagonal matrix,
$$\Lambda_i \mathop  = \limits^\Delta diag \left( {{( {\sigma _{{{\left\langle {i} \right\rangle }_t}}^{\left( t \right)}} )^2},{( {\sigma _{{{\left\langle m+1 \right\rangle }_t}}^{\left( t \right)}})^2},\cdots,{( {\sigma _{{{\left\langle k \right\rangle }_t}}^{\left( t \right)}})^2}}\right)~.$$

\begin{figure}[htbp]
\center{\includegraphics[width=0.43\textwidth]{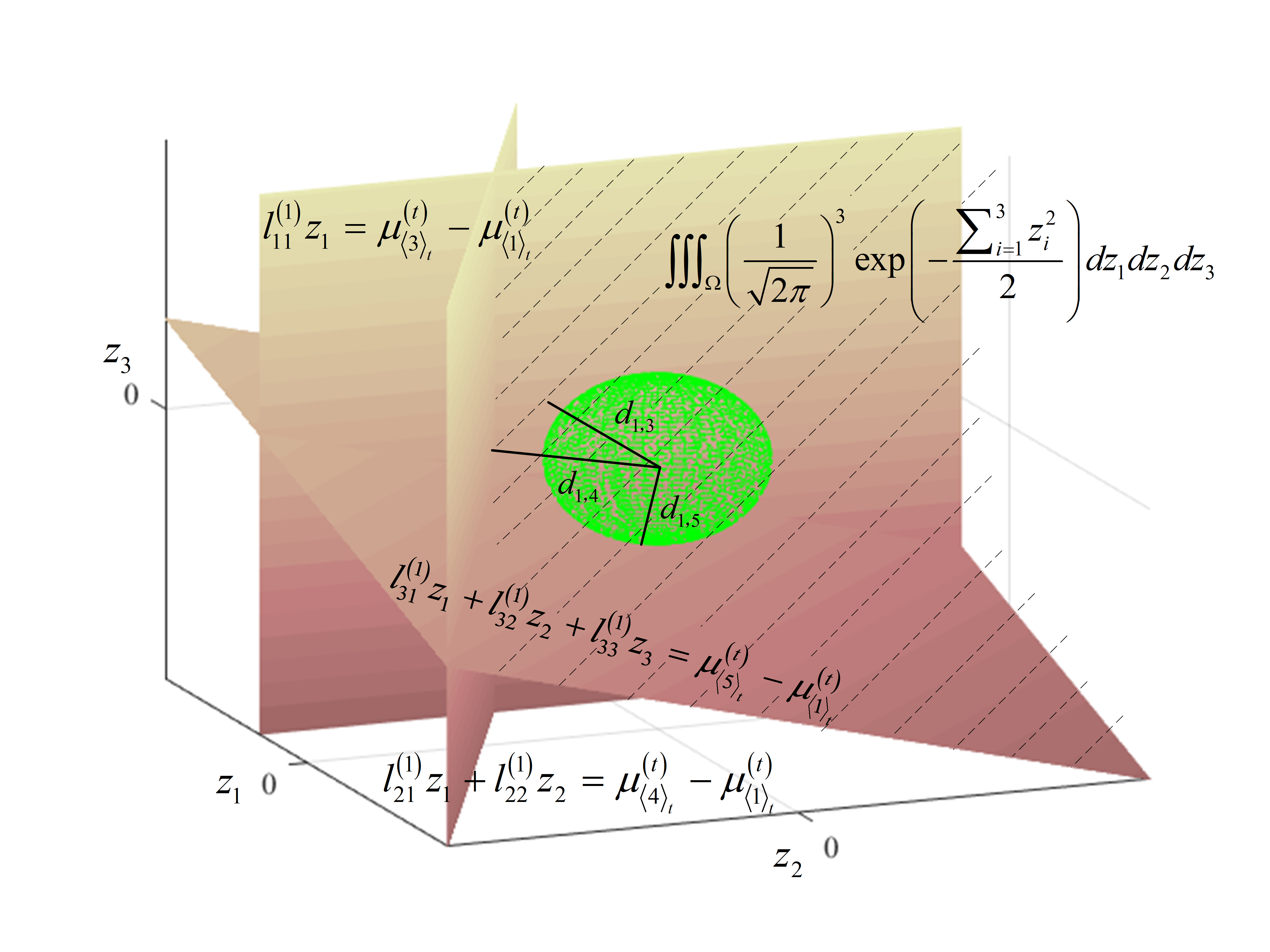}}
\caption{Area of integration for approximation is the ball, where dominant values of integrand $\exp\left(-\left(z_1^2+z_2^2+z_3^2\right)/2\right)$ are captured.}
\label{fig1}
\end{figure}

In addition, for $i \in \left\{1,\cdots,m\right\}$,
\begin{align}\label{inthyper}
& \Pr \left\{ {\left. {{\mu _{{{\left\langle {i} \right\rangle }_t}}} > {\mu _{{{\left\langle j \right\rangle }_t}}},\;j = m + 1, \cdots ,k} \right|{\mathcal{E} _t}} \right\} \\
= & \Pr \left\{\sum\limits_{h = 1}^{j - m} {{\ell_{\left( {j - m} \right),h}^{(i)}}{z_h}}> \mu _{{{\left\langle j \right\rangle }_t}}^{\left( t \right)} - \mu _{{{\left\langle {i} \right\rangle }_t}}^{\left( t \right)}, j = m + 1, \cdots ,k \right\} \nonumber\\
= & \frac{1}{{{{\left( {2\pi } \right)}^{{{\left( {k - m} \right)} \mathord{\left/{\vphantom {{\left( {k - 1} \right)} 2}} \right.\kern-\nulldelimiterspace} 2}}}}} \times \int_{\sum\nolimits_{h = 1}^{j - m} {{\ell_{\left( {j - m} \right),h}^{(i)}}} {z_h} > \mu _{{{\left\langle j \right\rangle }_t}}^{\left( t \right)} - \mu _{{{\left\langle {i} \right\rangle }_t}}^{\left( t \right)}} \nonumber\\
 & {\exp\left( {\frac{{-\sum\nolimits_{{\widetilde h} = 1}^{k - m} {z_{\widetilde h}^2} }}{2}}\right)d{z_1} \cdots d{z_{ {k - m} }},\;j=m+1,\cdots,k}~, \nonumber
\end{align}
where $L_i \mathop  = \limits^\Delta [ {{\ell_{\left( {j - m} \right),h}^{(i)}}} ] \in {\mathbb{R}^{\left( {k - m} \right) \times \left( {k - m} \right)}}$, $j=m+1,\cdots,k$, $h=1,\cdots,k-m$ is a lower triangular matrix of Cholesky decomposition for $\Sigma_i  = {L_{i}}L_{i}'$, and ${z_h}\mathop \sim \limits^{i.i.d} N\left( {0,1} \right)$. Therefore, (\ref{inthyper}) is an integral of the density of $\left(k-m\right)$ dimensional standard normal distribution over a region formed by some hyperplanes, and the value function becomes the maximum value of $m$ integrals. We visualize (\ref{inthyper}) and its approximation when $k=5$, $m=2$ and $i=1$. In Figure \ref{fig1}, (\ref{inthyper}) is the integration of the three-dimensional standard normal density over the shadowed area, and an approximation is captured by the size of the largest internally tangent ball with the radius $d_{1,5}$. In general, the VFA is given by ${\widetilde {V}_t}\left( {{\mathcal{E} _t}} \right) \mathop {\rm{ = }}\limits^\Delta \mathop {\max }\nolimits_{i = 1, \cdots ,m} d_{i}^2\left( {{\mathcal{E} _t}} \right)$, such that ${d_{i}}\left( {{\mathcal{E} _t}} \right) = \mathop {\min }  \left( {{d_{{i},\left(m + 1\right)}}\left( {{\mathcal{E} _t}} \right), \cdots ,{d_{{i},k}}\left( {{\mathcal{E} _t}} \right)} \right)$, where
$${d_{{i},j}}\left( {{\mathcal{E} _t}} \right) = \frac{{\mu _{{{\left\langle {i} \right\rangle }_t}}^{\left( t \right)} - \mu _{{{\left\langle j \right\rangle }_t}}^{\left( t \right)}}}{{\sqrt {{\left( {\sigma _{{{\left\langle {i} \right\rangle }_t}}^{\left( t \right)}} \right)^2} + {\left( {\sigma _{{{\left\langle j \right\rangle }_t}}^{\left( t \right)}} \right)^2}} }},\;\tabincell{c}{$i=1,\cdots,m$ \\ $j=m+1,\cdots,k$}~.$$

As $t_h \to \infty$, by the law of large numbers (LLN), $\mathop {\lim }\nolimits_{t \to \infty } \mu _{h}^{\left(t\right)}  = {\mu _h},\;a.s.$ and $\mathop {\lim }\nolimits_{t \to \infty } \sigma _{h}^{\left(t\right)} = 0,\;a.s.$. If alternatives $i$ and $j$ are sampled infinitely often as $t$ goes to infinity, $\mathop {\lim }\nolimits_{t \to \infty } {d_{i,j}}\left( {{\mathcal{E} _t}} \right) = \infty$. The VFA is reasonable since the exponential decreasing rate of the normal density.

\begin{prop}\label{apperr0}
 The error of the integral of the $\left(k-m\right)$ dimensional standard normal density over a centered ball with a radius $d_i\left(\mathcal{E}_t\right)$ as an approximation of (\ref{inthyper}) decreases to 0 in an exponential rate as $d_i\left(\mathcal{E}_t\right) \to \infty$, $i \in \left\{1,\cdots,m\right\}$.
\end{prop}

{The proof of the proposition can be found in the online appendix \cite{zhang2022online}.} If the $\left(t+1\right)$-th replication is the last one, a VFA one-step look ahead is given as ${\widetilde {V}_t}\left( {{\mathcal{E} _t};i} \right)\mathop  = \limits^\Delta  \mathbb{E} [ { {{{\widetilde {V}_{t+1}}}\left( {{\mathcal{E} _t} \cup \left\{ {{X_{i,t_{i} + 1}}} \right\}} \right)} |{\mathcal{E} _t}} ]$. With a certainty equivalence which replaces stochastic quantities by their expected values~\cite{Bertsekas2005}, a further approximation to soften the computational difficulty of the above expectation can be given as ${\widehat {V}_t}\left( {{\mathcal{E} _t};i} \right)\mathop  = \limits^\Delta {\widetilde {V}_{t+1}}\left( {{\mathcal{E} _t} \cup \mathbb{E}\left[ {\left. {{X_{i,t_{i} + 1}}} \right|{\mathcal{E} _t}} \right]} \right)$, $i=1,\cdots,k$. We use ${\widehat {V}_t}\left( {{\mathcal{E} _t};i} \right)$ as the criterion for making dynamic sampling decisions. For $i, h = 1,\cdots,m$, $j, \ell = m+1,\cdots,k$, and $i \ne h$, $j \ne \ell$,
$$\begin{array}{l}
{\widehat {V}_t}\left( {{\mathcal{E} _t};i} \right)=\\
\mathop {\max}\limits_{i,h} \left( {\mathop {\min }\limits_j \frac{{{{\left( {\mu _{{{\left\langle {i} \right\rangle }_t}}^{\left( t \right)} - \mu _{{{\left\langle j \right\rangle }_t}}^{\left( t \right)}} \right)}^2}}}{{{{\left( {\sigma _{{{\left\langle {i} \right\rangle }_t}}^{\left( {t + 1} \right)}} \right)}^2} + {{\left( {\sigma _{{{\left\langle j \right\rangle }_t}}^{\left( t \right)}} \right)}^2}}},\mathop {\min }\limits_{j} \frac{{{{\left( {\mu _{{{\left\langle {h} \right\rangle }_t}}^{\left( t \right)} - \mu _{{{\left\langle j \right\rangle }_t}}^{\left( t \right)}} \right)}^2}}}{{{{\left( {\sigma _{{{\left\langle {h} \right\rangle }_t}}^{\left( t \right)}} \right)}^2} + {{\left( {\sigma _{{{\left\langle j \right\rangle }_t}}^{\left( t \right)}} \right)}^2}}}} \right)~,
\end{array}$$
$$\begin{array}{l}
{\widehat {V}_t}\left( {{\mathcal{E} _t};j} \right)=\\
\mathop{\max}\limits_{i} \left( {\mathop {\min }\limits_{j,\ell} \left( {\frac{{{{\left( {\mu _{{{\left\langle {i} \right\rangle }_t}}^{\left( t \right)} - \mu _{{{\left\langle j \right\rangle }_t}}^{\left( t \right)}} \right)}^2}}}{{{{\left( {\sigma _{{{\left\langle {i} \right\rangle }_t}}^{\left( t \right)}} \right)}^2} + {{\left( {\sigma _{{{\left\langle j \right\rangle }_t}}^{\left( {t + 1} \right)}} \right)}^2}}},\mathop {\min }\limits_{\ell} \frac{{{{\left( {\mu _{{{\left\langle {i} \right\rangle }_t}}^{\left( t \right)} - \mu _{{{\left\langle \ell \right\rangle }_t}}^{\left( t \right)}} \right)}^2}}}{{{{\left( {\sigma _{{{\left\langle {i} \right\rangle }_t}}^{\left( t \right)}} \right)}^2} + {{\left( {\sigma _{{{\left\langle \ell \right\rangle }_t}}^{\left( t \right)}} \right)}^2}}}} \right)} \right),
\end{array}$$
and
$$( \sigma _{i}^{\left( t+1 \right)} )^2 = {\left( {\frac{1}{{{( {\sigma _{i}^{( 0 )}} )^2}}} + \frac{{{t_{i}+1}}}{\sigma_{i}^2}} \right)^{ - 1}},\;i=1,\cdots,k~.$$

A fully sequential allocation procedure is given by
\begin{equation}\label{AOAge3}
{\widehat{A}_{t + 1}}\left( {{\mathcal{E} _t}} \right) {{\in}} \arg \mathop {\max }\nolimits_{i = {1,\cdots,k}} \widehat {{V_t}}\left( {{\mathcal{E} _t};i} \right)~.
\end{equation}

The VFA can be rewritten as $\mathop {\max }\limits_{i = 1, \cdots,m} \mathop {\min }\limits_{j = m + 1, \cdots ,k} {1 \mathord{\left/{\vphantom {1 {c_v^2\left( {i,j} \right)}}} \right.\kern-\nulldelimiterspace} {c_v^2\left( {i,j} \right)}}$, where ${c_v}\left( {i,j} \right)$ is the posterior noise-signal ratio (coefficient of variation) of ${\mu _{{{\left\langle {i} \right\rangle }_t}}} - {\mu _{{{\left\langle j \right\rangle }_t}}}$. The larger is ${c_v}\left( {i,j} \right)$,  the more difficult to correctly identify the sign of the difference of ${\mu _{{{\left\langle {i} \right\rangle }_t}}} - {\mu _{{{\left\langle j \right\rangle }_t}}}$. (\ref{AOAge3}) sequentially allocates each replication to an alternative to minimize the largest noise-signal ratio for pairs of alternatives in comparison for selecting a good enough subset. We then analyze the asymptotic property of (\ref{AOAge3}). {Let $\widehat {\mathcal F}_t^{k - m} \mathop  = \limits^\Delta \left\{ {{{\left\langle {m + 1} \right\rangle }_t}, \cdots ,{{\left\langle k \right\rangle }_t}} \right\}$. Define $\mathcal{T} \mathop  = \limits^\Delta \{ {i:i \in \bigcap\nolimits_{n = 1}^\infty  {\bigcup\nolimits_{t = n}^\infty  {\widehat {\mathcal{F}}_t^m} } } \}$ containing alternatives that appear infinitely often among top-$m$ estimated set, and $\mathcal{B} \mathop  = \limits^\Delta \{ {i:i \in \bigcap\nolimits_{n = 1}^\infty  {\bigcup\nolimits_{t = n}^\infty  {{\widehat {\mathcal F}}_t^{k - m}} } } \}$ containing alternatives that appear infinitely often among $k$ alternatives outside of the top-$m$ estimated set.} The proposed sampling rule is proved to be consistent in the following theorem.

\begin{thm}\label{consistent}
As $t \to \infty$, the sampling rule (\ref{AOAge3}) is consistent, i.e., $\mathop {\lim }\nolimits_{t \to \infty } {\left\langle 1 \right\rangle } \in {{\widehat {\mathcal{F}}}_t^m},\;a.s.$
\end{thm}

\begin{IEEEproof}
{We show that following the sampling rule (\ref{AOAge3}), at least an alternative in $\mathcal{T}$ and all of the $\left(k-m\right)$ alternatives in $\mathcal{B}$ will be sampled infinitely often almost surely}. Then by the LLN, at least one alternative performs better than $\left(k-m\right)$ alternatives, and the consistency holds.

Define $\Phi \mathop  = \limits^\Delta  \{ i: \mbox{alternative $i$ is sampled infinitely often, }a.s. \}$. If $\mathcal{T} \cap \Phi=\emptyset$ and $\mathcal{B} \cap \Phi \neq \emptyset$, then $\forall~i\in\mathcal{T}$,
$$\mathop {\lim }\limits_{t \to \infty } {( {\sigma _i^{\left( t \right)}})^2} > 0,\quad \mathop {\lim }\limits_{t \to \infty } [{{( {\sigma _{i}^{\left( t \right)}})^2} - ( {\sigma _{i}^{\left( t+1 \right)}})^2}]> 0~,$$
and $\mathop {\lim}\nolimits_{t \to \infty } [{{\widehat{V}_t} \left( {{\mathcal{E} _t};i} \right) - {\widetilde {V}_t}\left( {{\mathcal{E} _t}} \right)} ] > 0,\;a.s.$; $\exists~ j \in \mathcal{B}$ such that
$$\mathop {\lim }\limits_{t \to \infty } {( {\sigma _j^{\left( t \right)}})^2} =0, \quad \mathop {\lim }\limits_{t \to \infty } [{{( {\sigma _{j}^{\left( t \right)}})^2} - ( {\sigma _{j}^{\left( t+1 \right)}})^2}]=0~,$$
and $\mathop {\lim}\nolimits_{t \to \infty } [{{\widehat{V}_t} \left( {{\mathcal{E} _t};j} \right) - {\widetilde {V}_t}\left( {{\mathcal{E} _t}} \right)} ] =0,\;a.s.$, {which leads to a contradiction to the sampling rule (\ref{AOAge3}) that the alternative with the largest VFA is sampled.} Therefore, $\mathcal{T} \cap \Phi \neq \emptyset$.

If $\mathcal{B} \cap \Phi=\emptyset$, then $\forall~j \in \mathcal{B}$, $\mathop {\lim}\nolimits_{t \to \infty } [{{\widehat{V}_t} \left( {{\mathcal{E} _t};j} \right) - {\widetilde {V}_t}\left( {{\mathcal{E} _t}} \right)} ] >0,\;a.s.$; $\exists~i\in \Phi\cap\mathcal{T}$, $\mathop {\lim}\nolimits_{t \to \infty } [{{\widehat{V}_t} \left( {{\mathcal{E} _t};i} \right) - {\widetilde {V}_t}\left( {{\mathcal{E} _t}} \right)} ] =0,\;a.s.$,
which leads to an contradiction to the sampling rule (\ref{AOAge3}). Therefore, $\mathcal{B}\cap\Phi \neq \emptyset$. In addition, if $\mathcal{B}\setminus\Phi\neq \emptyset$, where $\backslash$ denotes the set minus, then $\exists~i \in \mathcal{T}\cap\Phi$, $\exists~j \in \mathcal{B}\cap\Phi$, and $\exists~\ell \in \mathcal{B}\setminus\Phi$, such that $\lim\nolimits_{t\to \infty}d_{i,j}(\mathcal{E}_t)=\infty$, $\lim\nolimits_{t\to \infty}d_{i,\ell}(\mathcal{E}_t)>0$, and $\mathop {\lim}\nolimits_{t \to \infty } [{{\widehat{V}_t} \left( {{\mathcal{E} _t};j} \right) - {\widetilde {V}_t}\left( {{\mathcal{E} _t}} \right)} ] =0,\;a.s.$, $\mathop {\lim}\nolimits_{t \to \infty } [{{\widehat{V}_t} \left( {{\mathcal{E} _t};\ell} \right) - {\widetilde {V}_t}\left( {{\mathcal{E} _t}} \right)} ] >0,\;a.s.$, which leads to an contradiction to the sampling rule (\ref{AOAge3}). Therefore, $\mathcal{B}\subset\Phi$.

{Then we claim that $\left| \mathcal{B} \right| = k-m$, where $\left| \cdot \right|$ is the cardinality of a set. If $\left| \mathcal{B} \right| < k-m$, then $\forall~i \in \mathcal{B}^c$, where $\mathcal{B}^c = \{{i:i \in \bigcup\nolimits_{n = 1}^\infty  {\bigcap\nolimits_{t = n}^\infty  {{\widehat {\mathcal F}}_t^{m}} } } \}$ is the complement of the set $\mathcal{B}$ by DeMorgan's law, $\exists~n_i \ge 1$ such that $\forall~t \ge n_i$, $i \notin {\widehat {\mathcal F}_t^{k-m}}$. Let $\widetilde n = \mathop {\max }\nolimits_{i \in {{\mathcal B}^c}} {n_i}$, and then $\forall~t \ge {\widetilde n}$, $\forall~i \in \mathcal{B}^c$, $i \notin {\widehat {\mathcal F}_t^{k-m}}$, leading to ${\widehat {\mathcal F}_t^{k-m}} \subseteq \mathcal{B}$, which contradicts $| {\widehat {\mathcal F}_t^{k-m}} | = k-m > \left| \mathcal {B} \right|$. Thus, $\left| \mathcal {B} \right| \ge k-m$. If $\left| \mathcal{B} \right| > k-m$, with $\mathcal{B}\subset\Phi$, we have $\forall~j\in\mathcal{B}$, $\mathop {\lim }\nolimits_{t \to \infty } \mu _j^{\left( t \right)} = {\mu _j},\;a.s.$ Then $\forall~\varepsilon > 0$, $\exists~t_j > 0$ such that $\forall~t > t_j$, $| {\mu _j^{\left( t \right)} - {\mu _j}} | < \varepsilon,\;a.s.$ Let ${\widetilde t} = \mathop {\max }\nolimits_{j \in \mathcal{B}} {t_j}$, and for $j,\ell \in \mathcal{B}$, $j \ne \ell$, if the inequality $\mu_j > \mu_{\ell}$ holds, then $\forall~t > {\widetilde t}$ such that $\mu _j^{\left( t \right)} > \mu _{\ell}^{\left( t \right)}$, due to $\mu _j^{\left( t \right)} - \mu _\ell^{\left( t \right)} > {\mu _j} - \varepsilon  - \left( {{\mu _\ell} + \varepsilon } \right) > 0$ for small enough $\varepsilon$. Since $\left| \mathcal{B} \right| > k-m$, $\exists~j_0 \in \mathcal{B},~J \subset \mathcal{B}$ with $\left| J \right| = k-m$, such that the inequality ${\mu _{j_0}} >  {\mu_j}$ holds for $j \in J$, and then $\forall t >{\widetilde t}$, ${\mu _{j_0}^{\left( t \right)}} >  {\mu_j^{\left( t \right)}}$ for $j \in J$, which yields $\forall t > {\widetilde t}$, $j_0 \notin {\widehat {\mathcal F}_t^{k-m}}$, leading to a contradiction to $j_0 \in \mathcal{B}$. Therefore, $\left| \mathcal{B} \right| = k-m$.} Summarizing the above, the theorem is proved.
\end{IEEEproof}

Let ${G_{ij}} ( {{r_{i}},{r_{j}}} )\mathop  = \limits^\Delta  \frac{{{{\left( {{\mu _i} - {\mu _j}} \right)}^2}}}{{{{\sigma _i^2} \mathord{\left/{\vphantom {{\sigma _i^2} {{r_{i}}}}} \right.\kern-\nulldelimiterspace} {{r_{i}}}} + {{\sigma _j^2} \mathord{\left/{\vphantom {{\sigma _j^2} {{r_{j}}}}} \right.\kern-\nulldelimiterspace} {{r_{j}}}}}}$, $i \in \mathcal{T}$, $j \in \mathcal{B}$, where ${r_h}\mathop  = \limits^\Delta  {{{t_h}} \mathord{\left/{\vphantom {{{t_i}} t}} \right.
 \kern-\nulldelimiterspace} t}$, $h=1,\cdots,k$, is the sampling ratio. Denote $J_i\mathop  = \limits^\Delta \{\ell: {G_{i \ell}}\left( r_i,r_\ell \right)=\mathop {\min }\nolimits_{j \in \mathcal{B}} {G_{ij}}\left( r_{i},r_{j} \right) \}$. The following theorem further establishes that only $\left(k-m+1\right)$ alternatives will be sampled infinitely often almost surely.

\begin{thm}\label{only}
For any $\Gamma \subset \mathcal{T}$ and $|\Gamma|>1$, suppose conditions (\ref{overallequation1}), (\ref{overallequation2}) and (\ref{ratiobal}) cannot hold simultaneously. Then $\left| \Phi \right|=k-m+1$,
\begin{equation}\label{overallequation1}
\mathop {\min }\limits_{j \in \mathcal{B}} {G_{ij}}\left( {{r_i},{r_j}} \right) = \mathop {\min }\limits_{h\in\Gamma} \mathop{\min}\limits_{j \in \mathcal{B}} {G_{hj}}\left( {{r_h},{r_j}} \right),
\end{equation}
where $\mathop{\sum}\nolimits_{i\in\Gamma} r_i+\mathop{\sum}\nolimits_{j\in\mathcal{B}}r_j = 1$,
\begin{align}\label{overallequation2}
\mathop{\bigcup}\nolimits_{i\in \Gamma} J_i=\mathcal{B}~,
\end{align}
\begin{equation}\label{ratiobal}
\begin{aligned}
& {\left. {\frac{\partial }{{\partial x}}\mathop {\min }\limits_{j \in \mathcal{B}} {G_{ij}}\left( {x,{r_{ j }}} \right)} \right|_{x = {r_{ i }}}}\\
= & {\left. {\mathop {\min }\limits_{h \in \Gamma } \frac{\partial }{{\partial x}}\mathop {\min }\limits_{j \in \mathcal{B}} {G_{hj}}\left( {x,{r_{ j }}} \right)} \right|_{x = {r_{ h }}}},\quad \forall~i\in\Gamma~.
\end{aligned}
\end{equation}
\end{thm}

\begin{IEEEproof}
Based on the Theorem~\ref{consistent}, the proof boils down to the question of verifying only one alternative in $\mathcal{T}$ will be sampled infinitely often almost surely following (\ref{AOAge3}).

Let $r_h^{\left( t \right)}\mathop  = \limits^\Delta  {{{t_h}} \mathord{\left/{\vphantom {{{t_h}} {t}}} \right.\kern-\nulldelimiterspace} {t}}$, $\sum\nolimits_{h = 1}^k {r_h^{\left( t \right)}}  = 1$, $h=1,\cdots,k$. By the LLN, $\mathop {\lim }\nolimits_{t \to \infty } \mu _{ {h} }^{\left( t \right)} \to {\mu _{ h }}$, $h \in \Phi$. Since the asymptotic sampling ratios will be determined by the increasing order of $d_{ij}\left(\mathcal{E}_t\right)$ with respect to $t$, we replace $\mu _{{ h }}^{\left( t \right)}$ and $( {\sigma _{{ h }}^{\left( t \right)}} )^2$ with $\mu _{{ h }}$ and ${{\sigma _{{ h }}^2} / {{t_h}}}$ in $d_{ij}\left(\mathcal{E}_t\right)$ for simplicity of analysis. Note that $0 \le r_{h}^{\left(t\right)} \le 1$, $h=1,\cdots,k$, and $( {r_1^{\left( t \right)}, \cdots ,r_k^{\left( t \right)}} )$ is a bounded sequence. By the Bolzano-Weierstrass theorem~\cite{rudin1964principles}, there exists a subsequence of $( {r_1^{\left( t \right)}, \cdots ,r_k^{\left( t \right)}} )$ converging to $\left( {{{\widetilde r}_1}, \cdots ,{{\widetilde r}_k}} \right)$, where $\sum\nolimits_{\ell = 1}^k {{{\widetilde r}_\ell}}  = 1$, ${{{\widetilde r}_\ell}} \ge 0$. Without loss of generality, we assume that $( {r_1^{\left( t \right)}, \cdots, r_k^{\left( t \right)}})$ converges to $\left( {{{\widetilde r}_1}, \cdots, {{\widetilde r}_k}} \right)$; otherwise, the following argument is made over a subsequence. Notice that
$$\begin{aligned}
&\mathop {\lim }\limits_{t \to \infty } \left[ {\frac{{{{\left( {{\mu _{ i}} - {\mu _{ j}}} \right)}^2}}}{{{{\sigma _{ i}^2} \mathord{\left/{\vphantom {{\sigma _{ i}^2} {\left( {{t_{ i }} + 1} \right)}}} \right.\kern-\nulldelimiterspace} {\left( {{t_{ i }} + 1} \right)}} + {{\sigma _{ j}^2} \mathord{\left/{\vphantom {{\sigma _{ j}^2} {{t_{ j }}}}} \right.\kern-\nulldelimiterspace} {{t_{ j }}}}}} - \frac{{{{\left( {{\mu _{ i}} - {\mu _{ j}}} \right)}^2}}}{{{{\sigma _{ i}^2} \mathord{\left/{\vphantom {{\sigma _{ i}^2} {{t_{ i }}}}} \right.\kern-\nulldelimiterspace} {{t_{ i }}}} + {{\sigma _{ j}^2} \mathord{\left/{\vphantom {{\sigma _{ j}^2} {{t_{ j }}}}} \right.\kern-\nulldelimiterspace} {{t_{ j }}}}}}} \right]\\
& = \mathop {\lim }\limits_{t \to \infty } t\left[ {{G_{ij}} ( {r_{ i }^{\left( t \right)} + {1 \mathord{\left/{\vphantom {1 t}} \right.
 \kern-\nulldelimiterspace} t},r_{ j }^{\left( t \right)}} ) - {G_{ij}} ( {r_{ i }^{\left( t \right)},r_{ j }^{\left( t \right)}} )} \right] \\
& = \mathop {\lim }\limits_{t \to \infty } {\left. {{{\partial {G_{ij}}( {x,r_{ j  }^{\left( t \right)}})} \mathord{\left/{\vphantom {{\partial {G_{ij}}\left( {x,r_{ j  }^{\left( t \right)}} \right)} {\partial x}}} \right.\kern-\nulldelimiterspace} {\partial x}}} \right|_{x = r_{ i  }^{\left( t \right)}}} \\
& = {\left( {\frac{{{\sigma _{ i}}}}{{\widetilde r}_{ i }}} \right)^2}\frac{{{{\left( {{\mu _{ i}} - {\mu _{ j}}} \right)}^2}}}{{{\left( {{{\sigma _{ i}^2} \mathord{\left/{\vphantom {{\sigma _{ i}^2} {{\widetilde r}_{ i }}}} \right.\kern-\nulldelimiterspace} {{\widetilde r}_{ i }}} + {{\sigma _{ j}^2} \mathord{\left/{\vphantom {{\sigma _{ j}^2} {{\widetilde r}_{ j }}}} \right.\kern-\nulldelimiterspace} {{\widetilde r}_{ j }}}} \right)^2}}},\;\tabincell{c}{$i \in \mathcal{T}  \cap \Phi$ \\ $j \in \mathcal{B}$}~,
\end{aligned}$$
$$\begin{aligned}
&\mathop {\lim }\limits_{t \to \infty } \left[ {\frac{{{{\left( {{\mu _{ i}} - {\mu _{ j}}} \right)}^2}}}{{{{\sigma _{ i}^2} \mathord{\left/{\vphantom {{\sigma _{ i}^2} { {t_{ i }}  }}} \right.\kern-\nulldelimiterspace} { {t_{ i }}  }} + {{\sigma _{ j}^2} \mathord{\left/{\vphantom {{\sigma _{ j}^2} {{t_{ j }}}}} \right.\kern-\nulldelimiterspace} {{(t_{ j }+1)}}}}} - \frac{{{{\left( {{\mu _{ i}} - {\mu _{ j}}} \right)}^2}}}{{{{\sigma _{ i}^2} \mathord{\left/{\vphantom {{\sigma _{ i}^2} {{t_{ i }}}}} \right.\kern-\nulldelimiterspace} {{t_{ i }}}} + {{\sigma _{ j}^2} \mathord{\left/{\vphantom {{\sigma _{ j}^2} {{t_{ j }}}}} \right.\kern-\nulldelimiterspace} {{t_{ j }}}}}}} \right]\\
& = \mathop {\lim }\limits_{t \to \infty } t\left[ {{G_{ij}} ( {r_{ i }^{\left( t \right)},r_{ j }^{\left( t \right)} + {1 \mathord{\left/
 {\vphantom {1 t}} \right.
 \kern-\nulldelimiterspace} t}} ) - {G_{ij}} ( {r_{ i }^{\left( t \right)},r_{ j }^{\left( t \right)}} )} \right] \\
& = \mathop {\lim }\limits_{t \to \infty } {\left. {{{\partial {G_{ij}}( {r_{ i  }^{\left( t \right)},x} )} \mathord{\left/{\vphantom {{\partial {G_{ij}}\left( {r_{ i  }^{\left( t \right)},x} \right)} {\partial x}}} \right.\kern-\nulldelimiterspace} {\partial x}}} \right|_{x = r_{ j  }^{\left( t \right)}}} \\
& = {\left( {\frac{{{\sigma _{ j}}}}{{\widetilde r}_{ j }}} \right)^2}\frac{{{{\left( {{\mu _{ i}} - {\mu _{ j}}} \right)}^2}}}{{{{\left( {{{\sigma _{ i}^2} \mathord{\left/{\vphantom {{\sigma _{ i}^2} {{\widetilde r}_{ i }}}} \right.\kern-\nulldelimiterspace} {{\widetilde r}_{ i }}} + {{\sigma _{ j}^2} \mathord{\left/{\vphantom {{\sigma _{ j}^2} {{\widetilde r}_{ j }}}} \right.\kern-\nulldelimiterspace} {{\widetilde r}_{ j }}}} \right)}^2}}},\;\tabincell{c}{$i \in \mathcal{T}  \cap \Phi$ \\ $j \in \mathcal{B}$}~.
\end{aligned}$$

We first claim that $\exists~i\in\mathcal{T}\cap\Phi$ such that $\widetilde{r}_{ i  }>0$; otherwise, $\exists$ $j\in\mathcal{B}$ such that  $\widetilde{r}_{ j  }>0$, and
$$\mathop {\lim }\limits_{t \to \infty } {\left. \frac{\partial {G_{ij}}( {x,r_{ j  }^{\left( t \right)}})}{\partial x} \right|_{x = r_{ i  }^{\left( t \right)}}} > 0, \; \mathop {\lim }\limits_{t \to \infty } {\left. \frac{\partial {G_{ij}}( {r_{ i  }^{\left( t \right)},x} )}{\partial x} \right|_{x = r_{ j  }^{\left( t \right)}}} = 0,$$
which contradicts the sampling rule (\ref{AOAge3}). We then claim that $\exists~j\in\mathcal{B}$ such that $\widetilde{r}_{j}>0$; otherwise, $\exists~i\in\mathcal{T}\cap\Phi$ such that  $\widetilde{r}_{i}>0$, and
$$\mathop {\lim }\limits_{t \to \infty } {\left. \frac{\partial {G_{ij}}( {x,r_{ j  }^{\left( t \right)}})}{\partial x} \right|_{x = r_{ i  }^{\left( t \right)}}} = 0, \; \mathop {\lim }\limits_{t \to \infty } {\left. \frac{\partial {G_{ij}}( {r_{ i  }^{\left( t \right)},x} )}{\partial x} \right|_{x = r_{ j }^{\left( t \right)}}} > 0,$$
which contradicts the sampling rule (\ref{AOAge3}). In addition, we claim that $\forall~j\in\mathcal{B}$ such that $\widetilde{r}_{ j  }>0$; otherwise, $\exists~\ell, q\in\mathcal{B}$, $\ell \ne q$ such that $\widetilde{r}_{ \ell  }>0$, $\widetilde{r}_{ q  }=0$, and for $i\in\mathcal{T}\cap\Phi$ such that ${\widetilde r_{ i  }} > 0$,
\begin{align*}
\mathop{\lim}\nolimits_{t\to\infty}{G_{i\ell}}( {{r_{ i }^{(t)}},{r_{ \ell }^{(t)}}} )>0,\; \mathop{\lim}\nolimits_{t\to\infty}{G_{iq}}( {{r_{ i }^{(t)}},{r_{ q }^{(t)}}} )=0~,
\end{align*}
which contradicts the sampling rule (\ref{AOAge3}). Finally, we claim that $\forall~i \in \mathcal{T} \cap \Phi$ such that ${\widetilde{r}_{ i  }} > 0$; otherwise, $\exists~h,p\in\mathcal{T} \cap \Phi$, $h \ne p$ such that $\widetilde{r}_{ h  } > 0$, $\widetilde{r}_{ p  } = 0$, and for $j \in \mathcal{B}$,
\begin{align*}
\mathop{\lim}\nolimits_{t\to\infty}{G_{hj}} ( {{r_{ h }^{(t)}},{r_{ j }^{(t)}}} )>0,\;\mathop{\lim}\nolimits_{t\to\infty}{G_{pj}} ( {{r_{ p }^{(t)}},{r_{ j }^{(t)}}} )=0~,
\end{align*}
which contradicts the sampling rule (\ref{AOAge3}). Summarizing the above, we have $\widetilde{r}_{ \ell  } >0$, $\ell \in \Phi$.

Following the sampling rule (\ref{AOAge3}), we will show that ${\widetilde r}_{ \ell }$, $\ell \in \Phi$ satisfy (\ref{overallequation1}) and (\ref{overallequation2}) with $\Gamma=\mathcal{T} \cap \Phi $.
If $|\mathcal{T}\cap \Phi|=1$, then (\ref{overallequation1}) holds automatically. If $|\mathcal{T}\cap\Phi|>1$ and (\ref{overallequation1}) does not satisfy, then $\exists~i,h \in \mathcal{T} \cap \Phi $, $i \ne h$, such that the inequality
$$\mathop {\min }\limits_{j \in \mathcal{B}} {G_{ij}}\left( {{{\widetilde r}_{ i }},{{\widetilde r}_{ j }}} \right) > \mathop {\min }\limits_{j \in \mathcal{B}} {G_{hj}}\left( {{{\widetilde r}_{ h }},{{\widetilde r}_{ j }}} \right)~,$$
strictly holds, and then there exists $T_0 > 0$, such that after step $t > T_0$,
$$\mathop {\min }\limits_{j \in\mathcal{B}}{G_{ij}} ( {{r_{ i }^{\left( t \right)}},{r_{ j }^{\left( t \right)}}} ) > \mathop {\min }\limits_{j \in \mathcal{B}} {G_{hj}} ( {{r_{ h }^{\left( t \right)}},r_{ j }^{\left( t \right)}} )~,$$
due to the continuity of the function $G_{ij}\left(r_i,r_j\right)$ with respect to $r_i,r_j,\mu_i,\mu_j,\sigma_i,\sigma_j$. By the sampling rule (\ref{AOAge3}), alternative $i$ will be sampled and alternative $h$ will stop receiving replications before the inequality reverses. Therefore, violation of condition (\ref{overallequation1}) leads to a contradiction to that $( {r_1^{\left( t \right)}, \cdots r_k^{\left( t \right)}})$ converges to $\left( {{{\widetilde r}_1}, \cdots {{\widetilde r}_k}} \right)$, which implies that (\ref{overallequation1}) must hold.

If (\ref{overallequation2}) does not hold, then $\exists~\ell\in\mathcal{B}$ such that for $i \in \mathcal{T}\cap\Phi$, the inequality
$$G_{i\ell}\left( {{{\widetilde r}_{ i }},{{\widetilde r}_{ \ell }}} \right)>\mathop{\min}\nolimits_{j\in \mathcal{B}}G_{ij}\left( {{{\widetilde r}_{ i }},{{\widetilde r}_{ j }}} \right)~,$$
strictly holds. By a similar argument,  alternative $ \ell $ will stop receiving replications before the inequality reverses, and then (\ref{overallequation2}) must hold.

If $|\mathcal{T}  \cap \Phi|>1$, then by the sampling rule (\ref{AOAge3}), ${\widetilde r}_{ \ell }$, $\ell \in \Phi$ satisfy (\ref{ratio1}) with $\Gamma=\mathcal{T} \cap \Phi $. Otherwise, $\exists~i,h\in \mathcal{T} \cap \Phi$ such that
$$\begin{aligned}
& \mathop {\lim }\limits_{t \to \infty } \left[ {\mathop {\min }\limits_{j \in \mathcal{B}} \frac{{{\left( {{\mu _{ i  }} - {\mu _{ j  }}} \right)^2}}}{{{{\sigma _{ i  }^2} \mathord{\left/{\vphantom {{\sigma _{ i  }^2} {\left( {{t_{ i }} + 1} \right)}}} \right.\kern-\nulldelimiterspace} {\left( {{t_{ i }} + 1} \right)}} + {{\sigma _{ j  }^2} \mathord{\left/{\vphantom {{\sigma _{ j  }^2} {{t_{ j }}}}} \right.\kern-\nulldelimiterspace} {{t_{ j }}}}}} - \mathop {\min }\limits_{j \in \mathcal{B}} \frac{{{\left( {{\mu _{ i  }} - {\mu _{ j  }}} \right)^2}}}{{{{\sigma _{ i  }^2} \mathord{\left/{\vphantom {{\sigma _{ i  }^2} {{t_{ i }}}}} \right.
 \kern-\nulldelimiterspace} {{t_{ i }}}} + {{\sigma _{ j  }^2} \mathord{\left/{\vphantom {{\sigma _{ j  }^2} {{t_{ j }}}}} \right.
 \kern-\nulldelimiterspace} {{t_{ j }}}}}}} \right] \\
 = & \mathop {\lim }\limits_{t \to \infty } {\left. {\frac{{\partial }}{{\partial x}}}\min_{j \in \mathcal{B}}{G_{ij}}( {x,r_{ j }^{\left( t \right)}}) \right|_{x = r_{ i }^{\left( t \right)}}}
\\
 >  & \mathop {\lim }\limits_{t \to \infty } {\left. {\frac{{\partial }}{{\partial x}}}\min_{j \in \mathcal{B}}{G_{hj}}( {x,r_{ j }^{\left( t \right)}}) \right|_{x = r_{ h }^{\left( t \right)}}} \\
 = & \mathop {\lim }\limits_{t \to \infty } \left[ {\mathop {\min }\limits_{j \in \mathcal{B}} \frac{{{{\left( {{\mu _{ h  }} - {\mu _{ j  }}} \right)}^2}}}{{{{\sigma _{ h  }^2} \mathord{\left/{\vphantom {{\sigma _{ h  }^2} {\left( {{t_{ h }} + 1} \right)}}} \right.\kern-\nulldelimiterspace} {\left( {{t_{ h }} + 1} \right)}} + {{\sigma _{ j  }^2} \mathord{\left/
 {\vphantom {{\sigma _{ j  }^2} {{t_{ j }}}}} \right.
 \kern-\nulldelimiterspace} {{t_{ j }}}}}} -  \mathop {\min }\limits_{j\in \mathcal{B}} \frac{{{{\left( {{\mu _{ h  }} - {\mu _{ j  }}} \right)}^2}}}{{{{\sigma _{ h  }^2} \mathord{\left/
 {\vphantom {{\sigma _{ h  }^2} {{t_{ h }}}}} \right.\kern-\nulldelimiterspace} {{t_{ h }}}} + {{\sigma _{ j  }^2} \mathord{\left/{\vphantom {{\sigma _{ j  }^2} {{t_{ j }}}}} \right.
 \kern-\nulldelimiterspace} {{t_{ j }}}}}}} \right]~,
\end{aligned}$$
which is contradictory to the sampling rule (\ref{AOAge3}). By assumption, (\ref{overallequation1}), (\ref{overallequation2}) and (\ref{ratiobal}) cannot hold simultaneously with $\Gamma=\mathcal{T}  \cap \Phi$ and $|\mathcal{T}  \cap \Phi|>1$. Therefore, $\left| \Phi \right|=k-m+1$, which proves the theorem.
\end{IEEEproof}

\begin{remark}
{(\ref{overallequation1}) and (\ref{overallequation2}) adjust the ratios of the simulation replications allocated to some pair-wise comparison alternatives. (\ref{ratiobal}) balances the incremental ratio of the simulation replications allocated to $i\in\Gamma$}. The condition of Theorem 2 can be checked straightforwardly in a special case when $\sigma_1=\cdots=\sigma_k$. In general, the condition could be checked numerically.
\end{remark}

\begin{coro}\label{coro}
Following the allocation procedure (\ref{AOAge3}), the sampling ratio of each alternative satisfy,
$$\mathop {\lim }\limits_{t \to \infty } r_i^{\left( t \right)} = \left\{
\begin{array}{ccc}
r_i & & {i \in \Phi} \\
0 & & \rm{else}
\end{array} \right. a.s.~,$$
where $\sum\nolimits_{i = 1}^k {r_i}  = 1$, and for $i' = \mathcal{T}  \cap \Phi$, $j,\ell \in \mathcal{B}$,
\begin{equation}\label{ratio1}
\frac{{{{\left( {{\mu _{ i' }} - {\mu _{ j }}} \right)}^2}}}{{{{\sigma _{ i' }^2} \mathord{\left/
 {\vphantom {{{{\left( {{\sigma _{ i' }}} \right)}^2}} {r_{ i' }}}} \right.
 \kern-\nulldelimiterspace} {r_{ i' }}} + {{\sigma _{ j }^2} \mathord{\left/
 {\vphantom {{{{\left( {{\sigma _{ j }}} \right)}^2}} {r_{ j }}}} \right.
 \kern-\nulldelimiterspace} {r_{ j }}}}} = \frac{{{{\left( {{\mu _{ i' }} - {\mu _{ \ell }}} \right)}^2}}}{{{{\sigma _{ i' }^2} \mathord{\left/
 {\vphantom {{{{\left( {{\sigma _{ i' }}} \right)}^2}} {r_{ i' }}}} \right.
 \kern-\nulldelimiterspace} {r_{ i' }}} + {{\sigma _{ \ell }^2} \mathord{\left/
 {\vphantom {{{{\left( {{\sigma _{ \ell }}} \right)}^2}} {r_{ \ell }}}} \right.
 \kern-\nulldelimiterspace} {r_{ \ell }}}}}~,
 \end{equation}
 and
\begin{equation}\label{ratio2}
 r_{ i' } = {{{\sigma _{ i' }}}}\sqrt {\sum\nolimits_{j \in \mathcal{B}} {{{{r_{ j }^2 }} \mathord{\left/{\vphantom {{{{\left( {r_{ j }} \right)}^2}} {\sigma _{ j }^2}}} \right.\kern-\nulldelimiterspace} {\sigma _{ j }^2}}} }~.
 \end{equation}
\end{coro}
\begin{remark}\label{rmk01}
The proof of Corollary \ref{coro} is essentially the same as the proof for the Theorem 3 in~\cite{peng2018ranking}, which establishes that the sampling ratios determined by (\ref{ratio1}) and (\ref{ratio2}) are asymptotically optimal for selecting the best alternative. {The proposed sequential sampling policy (\ref{AOAge3}) for selecting  a good enough subset is denoted as AOA-gs}. AOA-gs procedure dynamically determines the alternative that stops receiving simulation replications based on collected sample information. {We caution that the sets $\mathcal{T}$, $\mathcal{B}$, and the alternative $\mathcal{T} \cap \Phi$ cannot be uniquely determined in our theory, since they depend on the realization of simulation replications.} Compared with selecting the single best alternative, where all of the $k$ alternatives are sampled infinitely often, the sampling ratios of $\left(m-1\right)$ alternatives converge to 0 as $t \to \infty$, which significantly reduces the computational cost.
\end{remark}

\begin{prop}
There is only one solution to (\ref{ratio1}) and (\ref{ratio2}).
\end{prop}

\begin{IEEEproof}
Suppose ${{\widetilde r}^*} \mathop  = \limits^\Delta ( {\widetilde r_{ i'^{\mathcal{T}} }^*,\widetilde r_{ m+1^{\mathcal{B}} }^*,\widetilde r_{ m+2^{\mathcal{B}} }^*, \cdots ,\widetilde r_{ k^{\mathcal{B}} }^*} )$ and ${r^*} \mathop  = \limits^\Delta ( {r_{ i'^{\mathcal{T}} }^*,r_{ m+1^{\mathcal{B}} }^*,r_{ m+2^{\mathcal{B}} }^*, \cdots ,r_{ k^{\mathcal{B}} }^*} )$ are two different solutions to (\ref{ratio1}) and (\ref{ratio2}). Following \cite{glynn2004}, ${{\widetilde r}^*}$ and ${r^*}$ are derived by the Karush-Kuhn Tucker (KKT) conditions of a concave maximization optimization problem. Then $\lambda {{\widetilde r}^*} + \left( {1 - \lambda } \right){r^*}$, $\lambda \in \left[0,1\right]$ are solutions to (\ref{ratio1}) and (\ref{ratio2}) since the optimal set of a convex optimization problem is convex \cite{boyd2009convex}. By the chain rule, (\ref{ratio1}), (\ref{ratio2}) and $\sum\nolimits_{i \in \Phi } {r_i^*}  = 1$ lead to $\Psi \cdot \bar{R} = O$, where $\Psi \in \mathbb{R}^{\left(k-m+1\right)\times\left(k-m+1\right)}$ is defined in (\ref{equ16}), where for $i' = \mathcal{T}  \cap \Phi$, $j \in \mathcal{B}$,
$${\phi _{i'j}}\mathop  = \limits^\Delta  {\left. {\frac{{\partial {G_{i'j}} ( {r_{ i' }^*,y} )}}{{\partial y}}} \right|_{y = r_{ j }^*}},\quad {\psi _{i'j}}\mathop  = \limits^\Delta  {\left. {\frac{{\partial {G_{i'j}} ( {y,r_{ j }^*} )}}{{\partial y}}} \right|_{y = r_{ i' }^*}}~,$$
$${\psi _i}\mathop  = \limits^\Delta  \frac{{2r_{ i }^*}}{{\sigma _{ i }^2}},\quad {{\bar r}_{ i }}\mathop  = \limits^\Delta  r_{ i }^* - \widetilde r_{ i }^*,\quad i\in\Phi~,$$
$$\bar R\mathop  = \limits^\Delta  {\left( {{{\bar r}_{ m+1^{\mathcal{B}} }},{{\bar r}_{ m+2^{\mathcal{B}} }}, \cdots ,{{\bar r}_{ k^{\mathcal{B}} }},{{\bar r}_{ i'^{\mathcal{T}} }}} \right)^\prime } \in \mathbb{R}^{\left(k-m+1\right) \times {1}}~,$$
$$O\mathop  = \limits^\Delta  {\left( {0,0, \cdots ,0} \right)^\prime } \in \mathbb{R}^{\left(k-m+1\right) \times {1}}~.$$

Following \cite{glynn2004}, $\phi_{i'j} > 0$ and $\psi_{i'j} > 0$, $j \in \mathcal{B}$, and we have $\det \left( \Psi  \right) \ne 0$ and $\Psi$ is invertible, which yields $\bar{R} = O$. Then $r_i^* = {\widetilde r_i^*}$, $i \in \Phi$, which contradicts that $r^*$ and $\widetilde r^*$ are two different solutions to (\ref{ratio1}) and (\ref{ratio2}). Therefore, there is only one solution to (\ref{ratio1}) and (\ref{ratio2}), which proves the Proposition.
\newcounter{mytempeqncnt}
\begin{figure*}[htbp]
\normalsize
\setcounter{mytempeqncnt}{\value{equation}}
\setcounter{equation}{9}
\begin{equation}
\begin{aligned}\label{equ16}
\Psi \mathop  = \limits^\Delta  \left( {\begin{array}{*{20}{c}}
{{\phi _{i'\left( {m + 1} \right)}}}&{ - {\phi _{i'\left( {m + 2} \right)}}}& \cdots &0&0&{{\psi _{i'\left( {m + 1} \right)}} - {\psi _{i'\left( {m + 2} \right)}}}\\
0&{{\phi _{i'\left( {m + 2} \right)}}}& \cdots &0&0&{{\psi _{i'\left( {m + 2} \right)}} - {\psi _{i'\left( {m + 3} \right)}}}\\
 \vdots & \vdots & \ddots & \vdots & \vdots & \vdots \\
0&0& \cdots &{{\phi _{i'\left( {k - 1} \right)}}}&{ - {\phi _{i'k}}}&{{\psi _{i'\left( {k - 1} \right)}} - {\psi _{i'k}}}\\
{ - {\psi _{m + 1}}}&{ - {\psi _{m + 2}}}& \cdots &{ - {\psi _{k - 1}}}&{ - {\psi _k}}&{{\psi _{i'}}}\\
1&1& \cdots &1&1&1
\end{array}} \right)~.
\end{aligned}
\end{equation}
\setcounter{equation}{\value{mytempeqncnt}}
\hrulefill
\vspace*{4pt}
\end{figure*}
\end{IEEEproof}

\begin{remark}
{A closely related problem to soften the objective for finding the best is good alternative selection. Let $\mathcal{F}^m \mathop  = \limits^\Delta  \left\{ {\left\langle 1 \right\rangle , \cdots ,\left\langle m \right\rangle } \right\}$, and then the posterior PCS for finding a good enough alternative ${{{\left\langle {{i^*}} \right\rangle }_T}}$ can be expressed as
\begin{equation}\label{goodalter}
\begin{aligned}
& {\rm{PCS}}_T^{\prime} = \Pr \left\{ {\left. {{{\left\langle {{i^*}} \right\rangle }_T} \in {\mathcal{F}^m}} \right|{\mathcal{E}_T}} \right\} \\
& = \Pr \left\{ \left.{\bigcup\nolimits_{\scriptstyle J \subseteq \left\{ {1, \cdots ,k} \right\}\backslash {i^*}\atop
\scriptstyle \left| J \right| = k - m} {\bigcap\nolimits_{j \in J} {{\mu _{{{\left\langle {{i^*}} \right\rangle }_T}}} > {\mu _{{{\left\langle j \right\rangle }_T}}}} } }\right|{\mathcal{E}_T} \right\}~.
\end{aligned}
\end{equation}}

{The derivation of a dynamic allocation procedure for finding a good alternative based on (\ref{goodalter}) is essentially the same as the proposed sequential sampling rule (\ref{AOAge3}), however, the former needs to enumerate all of the possible cases of the set $J$ (in general, there are $C^{m-1}_{k-1}$ cases, where  $C^{m-1}_{k-1}$ denotes the number of combinations by choosing $\left(m-1\right)$ elements from $\left(k-1\right)$ elements), leading to a higher computational complexity than our method. Numerical examples to illustrate that for selecting a good alternative, the proposed allocation procedure (\ref{AOAge3}) has a comparable performance to the allocation procedure derived based on (\ref{goodalter}), can be found in the online appendix \cite{zhang2022online}.}
\end{remark}

\section{Numerical Experiments}\label{sec4}

In this section, we provide numerical experiments to demonstrate the efficiency of the proposed AOA-gs policy. {The detail of implementation for AOA-gs is shown in the online appendix \cite{zhang2022online}}. The equal allocation (EA), optimal computing budget allocation for $r$ good enough designs (OCBA-rgm) \cite{gao2017new}, and OCBA-rgm given the true parameter (OCBA-rgmt) are implemented for comparisons. EA equally allocates simulation replications to each alternative, i.e., roughly $T/ k$ replications for each alternative. OCBA-rgm derives simulation replications based on some asymptotic conditions, and $\left(m-r\right)$ alternatives stop receiving replications at each step. The details of implementation for OCBA-rgm can be found in online appendix \cite{zhang2022online}. OCBA-rgmt allocates simulation replications according to the OCBA-rgm rule under the perfect information (assuming the true parameters are known), whereas OCBA-rgm allocates simulation replications according to the OCBA-rgm rule with parameters estimated sequentially by available sample information. {OCBA-rgmt is compared to show that the difference between the asymptotic property and finite-sample performance of a sequential allocation procedure. The efficiency of each procedure is measured by the IPCS, which is reported as a function of sampling budget $T$, i.e., ${\rm{IPCS}}_T\mathop = \mathds{E}[ {{\mathds{1} ( {{{\left\langle {1} \right\rangle }} \in {{\widehat {\mathcal{F}}}_T^m}} )}} ]$}. In each numerical experiment, we set $n_0 = 10$ and IPCS is estimated by 100,000 independent macro experiments. The code for the numerical experiments in this note can be found in https://github.com/gongbozhang-pku/Good-Enough-Selection.

\textit{Experiment 1}. The allocation procedures are tested in a synthetic experiment with 50 competing designs. {The size of the good enough subset is $5$}. In each macro experiment, the performance of each alternative is generated from $\mu_i \sim N\left( {0,({{( {51 - i} )} \mathord{\left/{\vphantom {{\left( {51 - i} \right)} 2}} \right.\kern-\nulldelimiterspace} 10})^2} \right)$, $i=1,\cdots,50$. The simulation replications are drawn independently from a normal distribution $N\left({\mu_i},{\sigma_i^2}\right)$, where $\sigma_i = 51-i$. The total simulation budget is $T=1000$.

\begin{figure}[htbp]
\centering
\includegraphics[width=0.35\textwidth]{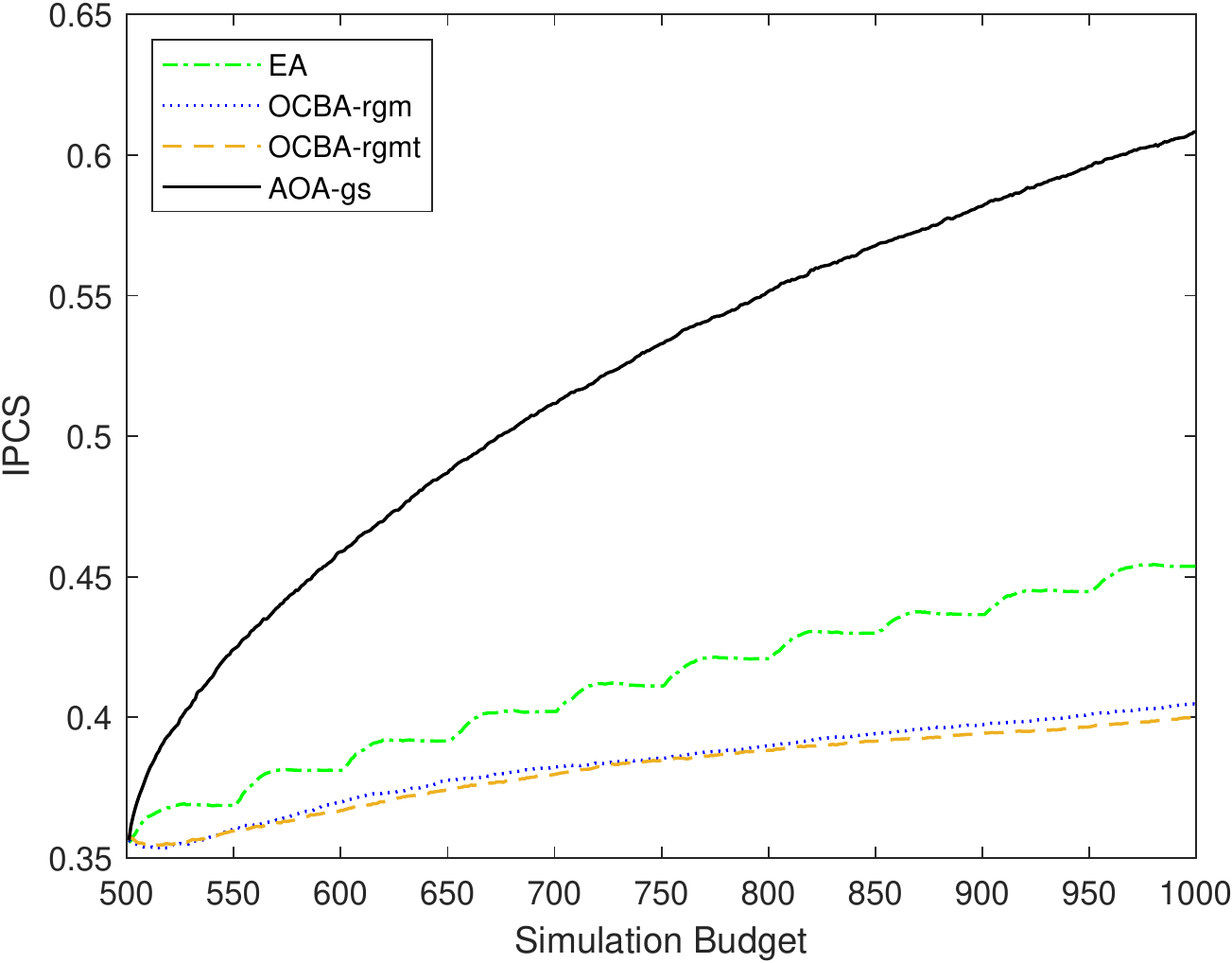}
\caption{Comparison of IPCS for four sampling allocation procedures in \textit{Experiment 1}.}
\label{fig3}
\end{figure}

{From Figure \ref{fig3}, we can see that OCBA-rgm has a slight edge over OCBA-rgmt that is worse than EA. AOA-gs achieves the best performance among all allocation procedures. The superior performance of AOA-gs could be attributed to the reason that AOA-gs makes full use of sample information and it is derived in a stochastic control framework with a support for finite sample performance.}

{We report the effects of $m$ on different allocation procedures in \textit{Experiment 2}. Table \ref{table1} reports the IPCS of four allocation procedures when $T=1000$ for $m=5,15,25,35,45$. From Table \ref{table1}, we can see that all allocation procedures obtain higher IPCS as $m$ grows, since the size of the good enough subset increases. EA performs better than OCBA-rgm when $m=5,15,25$. OCBA-rgm outperforms OCBA-rgmt, and AOA-gs stands out among all allocation procedures.}
\begin{table}[htbp]
\caption{IPCS of four allocation procedures for a variety of $m$ values}
\label{table1}
\centering
\begin{tabular}{cccccc}
\hline
   $m$ & EA & OCBA-rgm & OCBA-rgmt & AOA-gs  \\
\hline
   5  & 0.4538 & 0.4047 & 0.4000 & 0.6082\\
   15 & 0.7158 & 0.6072 & 0.6009 & 0.8481\\
   25 & 0.7964 & 0.7629 & 0.7451 & 0.9366 \\
   35 & 0.8633 & 0.9186 & 0.8713 & 0.9697\\
   45 & 0.9564 & 0.9837 & 0.9665 & 0.9915\\
\hline
\end{tabular}
\end{table}

{Three additional synthetic experiments, which highlight the difference between the asymptotic property and finite-sample performance for sequential OCBA-rgm procedure, and comparisons of AOA-gs and AOAP in \cite{peng2018ranking}, can be found in the online appendix \cite{zhang2022online}. Numerical results show that the advantage of AOA-gs appears to be more significant when the number of competing alternatives is large. }

\textit{Experiment 2: An $\left(s,S\right)$-type Inventory System Problem}. We consider an $\left(s,S\right)$-type inventory system problem. The goal is to determine the alternatives that have low expected average cost. The inventory system involves a single item under a periodic review with zero delivery lag and backlogging. The demand ${\xi _{t}}$ follows i.i.d. Poisson distribution with mean 25 and is independent of the level of inventory $\widetilde X_t$. The $\left(s,S\right)$-type inventory policy replenishes the stock at maximal capacity $S_i$ whenever the inventory level falls below a certain threshold $s_i$, that is,
$${\widetilde X_{t+1}} = \left\{
\begin{array}{lcl}
\widetilde X_t-{\xi _{t}}& , & {\xi _{t}}+s_i \le \widetilde X_t,\\
S_i  & ,& \rm{otherwise}.\\
\end{array} \right.$$

The total cost includes ordering cost (including a fixed setup cost $K$ of \$32 per order and an incremental cost $K^{\prime}$ of \$3 per stock), holding cost $H$ (\$1 per stock per stage) and shortage cost $p$ (\$5 per stock per stage). Parameters are set based on \cite{law2000simulation}. The cost at stage $t$ for $\left(s_i,S_i\right)$ policy is
$$\begin{aligned}
& C_i({\widetilde X}_t,\xi_t) = \\
& \left\{
\begin{array}{lcl}
H ({\widetilde X}_t-\xi_t )& , & {\xi _{t}} \le {\widetilde X}_t - s_i,\\
\left(H-K^{\prime}\right) ( {\widetilde X}_t-\xi_t ) + K + K^{\prime}S_i & ,&{\widetilde X}_t - s_i < {\xi _{t}} \le {\widetilde X}_t, \\
p (\xi_t-{\widetilde X}_t ) + K + K^{\prime}S_i&,& {\widetilde X}_t < \xi_t.
\end{array} \right.
\end{aligned}$$

20 alternatives are shown in Table \ref{table4}. The expected average cost is calculated over 30 stages. The best alternative is 12, i.e., $s_i = 20$, $S_i = 55$, determined by $10^7$ experiments for estimating the expected cost of each alternative. The size of the good enough subset is set as $3$. The efficiency of each procedure is measured by the PCS, since the alternatives have deterministic performances. The total simulation budget is $T=500$.

\begin{table}[htbp]
\caption{20 alternatives in an $\left(s,S\right)$-type inventory system}
\label{table4}
\centering
\begin{tabular}{ccccccccccc}
\hline
   Alternative & 1 & 2 & 3 & 4 & 5 & 6 & 7 & 8 & 9 & 10\\
\hline
   s & 5 & 5 & 10 & 10 & 10 & 10 & 10 & 10 & 20 & 20\\
   S & 45 & 50 & 45 & 50 & 55 & 60 & 65 & 70 & 40 & 45\\
\hline
   Alternative & 11 & 12 & 13 & 14 & 15 & 16 & 17 & 18 & 19 &20 \\
\hline
   s & 20 & 20 & 20 & 20 & 20 & 20 & 20 & 30 & 30 & 30\\
   S & 50 & 55 & 60 & 65 & 70 & 75 & 80 & 50 & 55 & 60\\
\hline
\end{tabular}
\end{table}

\begin{figure}[htbp]
\centering
\includegraphics[width=0.35\textwidth]{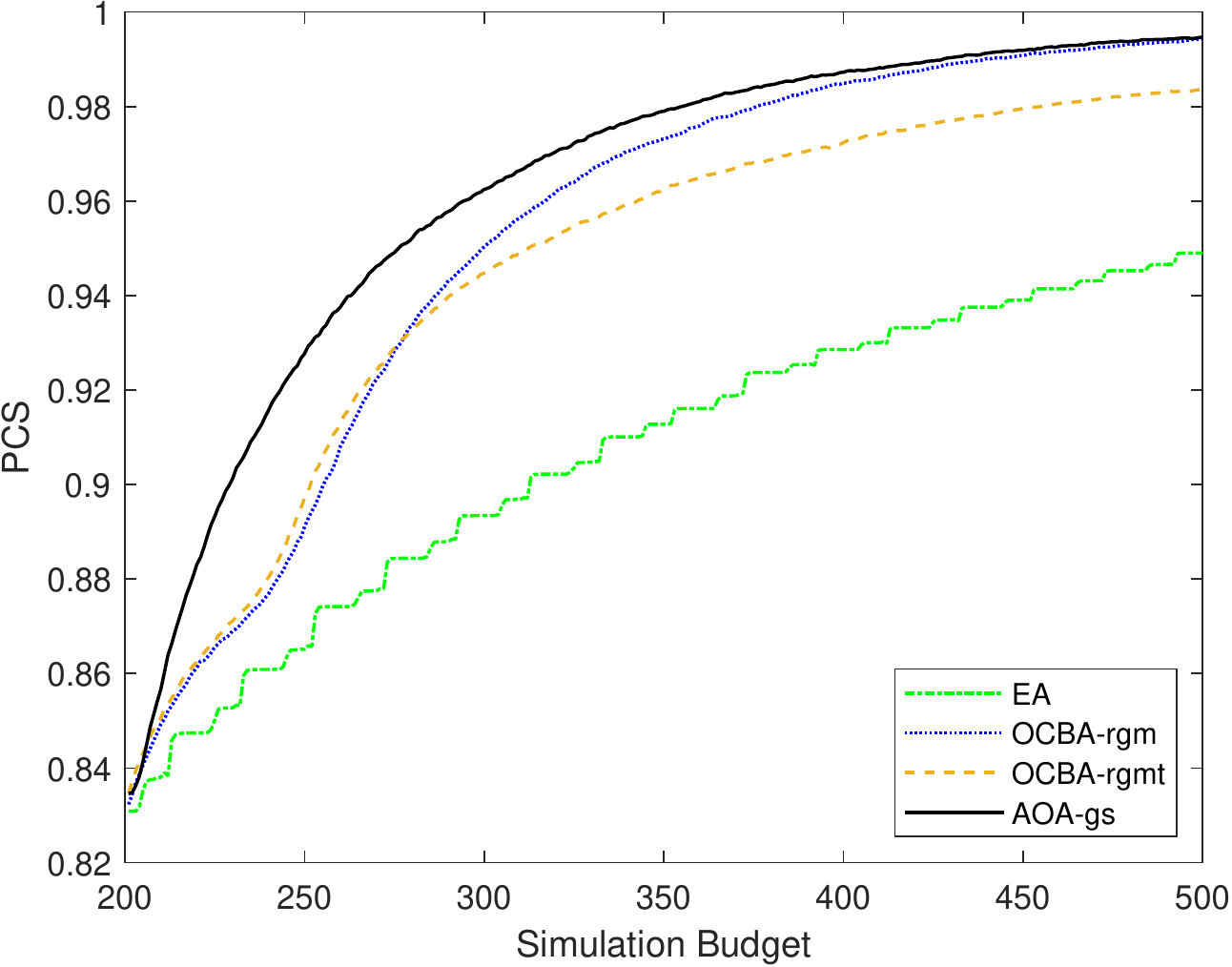}
\caption{Comparison of PCS for four allocation procedures in the $\left(s,S\right)$-type inventory system problem.}
\label{fig7}
\end{figure}

From Figure \ref{fig7}, we can see that EA performs the worst among all allocation procedures. OCBA-rgm and OCBA-rgmt have a comparable performance at the beginning, and the former increases at a faster rate than the latter as the simulation budget grows. AOA-ge achieves the best performance among all allocation procedures.

{An emergency department healthcare staffing problem can be found in the online appendix \cite{zhang2022online}.}

\section{Conclusion}\label{sec5}
{We present a sequential sampling procedure for selecting the good enough subset from a finite set of alternatives.  Given a simulation budget, the goal is to maximize the posterior PCS for selecting a good enough subset containing the best alternative.} We formulate the dynamic sampling decision under a stochastic control framework, and derive an efficient sequential allocation procedure named as AOA-gs by maximizing a VFA one-step look ahead.

The proposed sequential AOA-gs procedure is proved to be consistent, and we establish the asymptotic sampling ratios of AOA-gs, resulting in some useful insights. A set of experiments indicate that our proposed allocation procedure is significantly more efficient than existing sampling procedures, especially when the number of competing alternatives is large, under various sizes of the good enough subset. How to apply the proposed allocation procedure to improve the computational efficiency in reinforcement learning~\cite{jianglarge,jia2020structural,jiang2020computing} deserves future research.

\appendices

\section*{Acknowledgment}
This work was supported in part by the National Natural Science Foundation of China (NSFC) under Grants  72022001, 71901003, and 62073182. A preliminary study of this problem without detailed proofs has been presented in~\cite{zhang2020}.

\bibliographystyle{IEEEtran}
\bibliography{IEEEabrv,IEEEexample}

\end{document}